\documentclass[11]{article}
\usepackage[russian]{babel}
\usepackage{amssymb}
\usepackage{inputenc}
\usepackage{graphicx}
\begin{document}
\begin{center}
\textbf{\Large CALCULATION OF COEFFICIENTS OF THE OPTIMAL
QUADRATURE FORMULAS IN THE $W_2^{m,m-1}(0,1)$ SPACE}\\
\textbf{\large Kh.M.Shadimetov, A.R.Hayotov}
\end{center}

\textbf{Abstract.} In this paper problem of construction of
optimal quadrature formulas in $W_2^{(m,m-1)}(0,1)$ space is
considered. Here by using Sobolev's algorithm when $m=1,2$ we find
optimal coefficients of quadrature formulas of the form
$$
\int\limits_0^1\varphi(x)dx\cong
\sum\limits_{\beta=0}^NC_{\beta}\varphi(x_{\beta}).
$$

\large Many problems of science and engineering are reduced to a
integral and a differential equations or a system of such type
equations. Solutions of such equations is often expressed with the
help of defined integrals. But in many cases it is impossible to
calculate these integrals exactly. Therefore approximate
evaluation such integrals with possible high accuracy and less
expenditure is one of the \emph{actual problem} of computational
mathematics. This is well demonstrated in example of theory of
quadrature formulas. The task of finding of numerical value of one
dimensional integral, in view of its well-known geometrical
meaning, is often said the quadrature. There exist various methods
of quadrature, which allow to calculate defined integral with the
help of finite number values of integrand.

Present work also is devoted to one of such methods, i.e. to
construction of optimal quadrature formulas for approximate
calculation of defined integrals in the $W_2^{(m,m-1)}(0,1)$
space.

We consider following quadrature formula
$$
\int\limits_0^1p(x)\varphi(x)dx\cong \sum\limits_{\beta=0}^N
C_{\beta}\varphi(x_{\beta}) \eqno (1)
$$
with error functional
$$
\ell_N(x)=p(x)\varepsilon_{[0,1]}(x)- \sum\limits_{\beta=0}^N
C_{\beta}\delta(x-x_{\beta}) \eqno (2)
$$
on $W_2^{(m,m-1)}(0,1)$ space. Here $C_{\beta}$ are the
coefficients and $x_{\beta}$ are the nodes of the quadrature
formula (1), $\varepsilon_{[0,1]}(x)$ is the indicator of the
interval $[0,1]$, $p(x)$ is a weight function, $\delta(x)$ is
Dirac delta-function, $\varphi(x)$ is a element of the Hilbert
space $W_2^{(m,m-1)}(0,1)$, norm of a function in which is given
by formula
$$
\|\varphi(x)|W_2^{(m,m-1)}(0,1)\|=\left\{\int\limits_0^1\left(\varphi^{(m)}(x)+
\varphi^{(m-1)}(x)\right)^2dx\right\}^{1/2}.
$$
Difference
$$
(\ell_N(x),\varphi(x))= \int\limits_0^1p(x)\varphi(x)dx-
\sum\limits_{\beta=0}^N
C_{\beta}\varphi(x_{\beta})=\int\ell_N(x)\varphi(x)dx.
$$
is called by error of the quadrature formula (1).

Quality of the quadrature formula (1) is estimated by norm of the
error functional
$$
\|\ell_N(x)|W_2^{(m,m-1)*}(0,1)\|=\sup\limits_{\|\varphi(x)|W_2^{(m,m-1)}(0,1)\|=1}
|(\ell_N(x),\varphi(x))|.
$$
The norm of the error functional  $\ell_N(x)$ depends on the
coefficients $C_{\beta}$ and the nodes $x_{\beta}$. Choice of
coefficients in fixed nodes is linear problem. Therefore we will
fix the nodes  $x_{\beta}$ and the norm of the error functional
will minimize by the coefficients, i.e. we will find
$$
\|\stackrel{\circ}{\ell}_N(x)|W_2^{(m,m-1)*}(0,1)\|=\inf\limits_{C_{\beta}}
\|\ell_N(x)|W_2^{(m,m-1)*}(0,1)\|.\eqno(3)
$$

If it is found $\|\stackrel{\circ}{\ell}_N(x)\|$, then said, that
the functional $\stackrel{\circ}{\ell}_N(x)$ correspond to the
optimal quadrature formula in $W_2^{(m,m-1)}(0,1)$. Thus, we will
get following problems.

{\bf Problem 1.} {\it Find norm of the error functional
$\ell_N(x)$ of the quadrature formula of the form (1) in
$W_2^{(m,m-1)*}(0,1)$ space.}

{\bf Problem 2.} {\it Find such values of the coefficients
$C_{\beta}$, which satisfy the equality (3).}

Problem 1 was solved in [1] and for square of the norm of the
error functional  (2) was obtained
$$
\|\ell_N(x)\|^2=(-1)^m\Bigg[\sum\limits_{\beta=0}^N\sum\limits_{\beta'=0}^N
C_{\beta}C_{\beta'}\psi_m(x_{\beta}-x_{\beta'})-
$$
$$
-2\sum\limits_{\beta=0}^NC_{\beta}\int\limits_0^1p(x)\psi_m(x-x_{\beta})dx+
\int\limits_0^1\int\limits_0^1p(x)p(y)\psi_m(x-y)dxdy\Bigg],
\eqno(4)$$ where
$$
\psi_m(x)=\frac{\mathrm{sign}x}{2}\left(\frac{e^x-e^{-x}}{2}-
\sum\limits_{k=1}^{m-1}\frac{x^{2k-1}}{(2k-1)!}\right). \eqno(5)
$$
Furthermore, the error functional (2), as shown in [1], satisfies
following conditions
$$
(\ell_N(x),x^{\alpha})=0,\ \ \alpha=\overline{0,m-2}, \eqno(6)
$$
$$
(\ell_N(x),e^{-x})=0. \eqno(7)
$$

The norm  of the error functional is multidimensional function of
the coefficients $C_{\beta}$ $(\beta=\overline{0,N})$. For finding
point of condition minimum of the norm square of the error
functional (2) in conditions (6) and (7) in [1] applying the
method of Lagrange undetermined factors following system was
obtained
$$
\sum_{\gamma=0}^NC_{\gamma}\psi_m(x_{\beta}-x_{\gamma})+P_{m-2}(x_{\beta})+de^{-x_{\beta}}=\int\limits_0^1
p(x)\psi_m(x-x_{\beta})dx,\ \beta=\overline{0,N}, \eqno(8)
$$
$$
\sum_{\gamma=0}^N
C_{\gamma}x_{\gamma}^{\alpha}=\int\limits_0^1p(x)x^{\alpha}dx,\
\alpha=\overline{0,m-2}, \eqno(9)
$$
$$
\sum_{\gamma=0}^N
C_{\gamma}e^{-x_{\gamma}}=\int\limits_0^1p(x)e^{-x}dx, \eqno(10)
$$
where $\psi_m(x)$ is determined by formula (5), $P_{m-2}(x)$ is
unknown polynomial of degree $m-2$, $d$ is a constant.

In [1] also was proved existence and uniqueness of solution of the
system (8)-(10).

In present paper we will solve the system (8) - (10), i.e. we will
solve the problem 2.

Here interestingly that if system (8) - (10) is solved and the
optimal coefficients $C_{\beta}$ are found, then obtained
quadrature formula will be exact not only for polynomials of
degree $m-2$, it also will be exact for function $e^{-x}$.

Direct solution of the system (8)-(10) is very hard. Therefore we
will use Sobolev's algorithm [2], which is applied in solution
such systems.

Further, we briefly present this algorithm for the system (8) -
(10). For this we will determine $C_{\beta}$ for all integer
values of $\beta$, assume that for $\beta<0$ and $\beta>N$
coefficients $C_{\beta}$ are equal to zero. Let, furthermore,
$x_{\beta}=h\beta$, $h=\frac{1}{N}$, $N=1,2,....$

The system (8)-(10) by using the definition of convolution of
discrete functions we reduce to the following form
$$
\psi_m(h\beta)*C_{\beta}+P_{m-2}(h\beta)+de^{-h\beta}=f_m(h\beta),
\ h\beta\in [0,1], \eqno(11)
$$
$$
C_{\beta}=0 \mbox{ when } h\beta\not\in [0,1], \eqno(12)
$$
$$
\sum_{\beta=0}^NC_{\beta}(h\beta)^{\alpha}=q_{\alpha},\ \
\alpha=\overline{0,m-2}, \eqno(13)
$$
$$
\sum_{\beta=0}^NC_{\beta}e^{-h\beta}=q_{m-1},
\eqno(14)
$$
where $P_{m-2}(h\beta)$ is unknown polynomial of degree $m-2$, $d$
is a constant,
$$
f_m(h\beta)=\int\limits_0^1p(x)\psi_m(x-h\beta)dx,\eqno(15)
$$
$q_{\alpha}=\int\limits_0^1p(x) x^{\alpha}dx, $\
$(\alpha=\overline{0,m-2})$,
$q_{m-1}=\int\limits_0^1p(x)e^{-x}dx.$

Consider following problem.

\textbf{Problem A.} \emph{Find function $C_{\beta}$ and polynomial
 $P_{m-2}(h\beta)$, unknown coefficient $d$,
which satisfy the system (11) - (14) in given $f_m(h\beta)$ and}
$q_{\alpha}$ $(\alpha=\overline{0,m-1})$.

Idea of solution of the system (11) - (14), as in [2], consists on
replace unknown function $C_{\beta}$. Namely, instead of
$C_{\beta}$ we will introduce into consideration functions
$$
v(h\beta)=\psi_m(h\beta)*C_{\beta}$$ and
$$u(h\beta)=v(h\beta)+P_{m-2}(h\beta)+de^{-h\beta}.
$$
In such substitution we need only express $C_{\beta}$ with the
help of $u(h\beta)$. For this we must construct such operator
$D_m(h\beta)$, which satisfies following equality
$$
D_m(h\beta)*\psi_m(h\beta)=\delta(h\beta), \eqno(16)
$$
where $\delta(h\beta)$ is equal to zero when $h\beta\neq 0$ and is
equal to 1 when $h\beta=0$. In works [3,4] discrete operator
$D_m(h\beta)$ was constructed, its properties were investigated
and following theorems were proved

{\bf Theorem 1.} {\it Discrete analogue of the differential
operator $\frac{d^{2m}}{dx^{2m}}-\frac{d^{2m-2}}{dx^{2m-2}}$,
which satisfies equality (16) have following view:
$$
D_m(h\beta)=\frac{1}{p_{2m-2}^{(2m-2)}}\left\{
\begin{array}{ll}
\sum\limits_{k=1}^{m-1}A_k\lambda_k^{|\beta|-1}, & |\beta|\geq 2;\\
-2e^h+\sum\limits_{k=1}^{m-1}A_k, &|\beta|=1;\\
2C+\sum\limits_{k=1}^{m-1}\frac{A_k}{\lambda_k},& \beta=0,
\end{array}
\right. \eqno(17) $$ where
$$
C=1+(2m-2)e^h+e^{2h}+\frac{e^h\cdot
p_{2m-3}^{(2m-2)}}{p_{2m-2}^{(2m-2)}}, \eqno(18)
$$
$$
A_k=\frac{2(1-\lambda_k)^{2m-2}[\lambda_k(e^{2h}+1)-e^h(\lambda_k^2+1)]
p_{2m-2}^{(2m-2)}}{\lambda_kP_{2m-2}'(\lambda_k)}, \eqno(19)
$$
$$
P_{2m-2}(\lambda)=\sum\limits_{s=0}^{2m-2}p_s^{(2m-2)}\lambda^s=(1-e^{2h})
(1-\lambda)^{2m-2}-2(\lambda(e^{2h}+1)-e^h(\lambda^2+1))\times
$$
$$ \times
\left[h(1-\lambda)^{2m-4}+\frac{h^3(1-\lambda)^{2m-6}}{3!}E_2(\lambda)+
...+\frac{h^{2m-3}E_{2m-4}(\lambda)} {(2m-3)!}\right], \eqno(20)
$$
$p_{2m-2}^{(2m-2)}$, $p_{2m-3}^{(2m-2)}$ are coefficients of the
polynomial $P_{2m-2}(\lambda)$, $\lambda_k$ is a root of the
polynomial $P_{2m-2}(\lambda)$, $|\lambda_k|<1$, $E_k(\lambda)$ is
Euler polynomial [5].}

\textbf{Theorem 2.} \emph{The discrete analogue $D_m(h\beta)$ of
the differential operator
$\frac{d^{2m}}{dx^{2m}}-\frac{d^{2m-2}}{dx^{2m-2}}$ when $m=1,2,3$
satisfies following equalities}

1) $D_m(h\beta)*e^{h\beta}=0;$

2) $D_m(h\beta)*e^{-h\beta}=0;$

3) \emph{$D_m(h\beta)*(h\beta)^n=0,$ as $n\leq 2m-3$, i.e.
convolution of $D_m(h\beta)$ with polynomial of degree $\leq 2m-3$
when $m=2,3$ is equal to zero};

4) \emph{$D_m(h\beta)*\psi_m(h\beta)=\delta(h\beta),$ where
$\psi_m(h\beta)$ is determined by formula (5).}

Further, taking account of (16), we get
$$
C_{\beta}=D_m(h\beta)*u(h\beta). \eqno(21)
$$
So, if we find function $u(h\beta)$, then the optimal coefficients
are found from equality  (21).

Now we will find explicit form of the function $u(h\beta)$. From
equality (11) we get, that $u(h\beta)=f_m(h\beta)$ when $h\beta\in
[0,1]$. Since $C_{\beta}=0$ when $h\beta\not\in [0,1]$, then
$$
C_{\beta}=D_m(h\beta)*u(h\beta)=0\mbox{ when } h\beta\not\in
[0,1]. \eqno(22)
$$
We will calculate convolution $v(h\beta)=\psi_m(h\beta)*C_{\beta}$
when $h\beta\not\in [0,1]$.

Taking into account (12), (13), (14), we have $$ v(h\beta)=
\left\{
\begin{array}{ll}
-\frac{e^{h\beta}}{4}q_{m-1}+
De^{-h\beta}+Q^{(2m-3)}(h\beta)+Q_{m-2}(h\beta), & \beta<0;\\
\frac{e^{h\beta}}{4}q_{m-1}-
De^{-h\beta}-Q^{(2m-3)}(h\beta)-Q_{m-2}(h\beta), & \beta>N,
\end{array}
\right. \eqno(23)
$$
where
$$
Q^{(2m-3)}(h\beta)=
\frac{1}{2}\Bigg[\sum_{k=1}^{[\frac{m+1}{2}]-1}
\sum_{\alpha=0}^{2k-1}\frac{(h\beta)^{2k-1-\alpha}(-1)^{\alpha}}
{(2k-1-\alpha)!\cdot \alpha!}q_{\alpha}+
$$
$$ + \sum_{k=[\frac{m+1}{2}]}^{m-1}
\sum_{\alpha=0}^{m-2}\frac{(h\beta)^{2k-1-\alpha}(-1)^{\alpha}q_{\alpha}}
{(2k-1-\alpha)!\cdot \alpha!}\Bigg], \eqno(24) $$ $$
Q_{m-2}(h\beta)= \frac{1}{2}\sum_{k=[\frac{m+1}{2}]}^{m-1}
\sum_{\alpha=m-1}^{2k-1}\frac{(h\beta)^{2k-1-\alpha}(-1)^{\alpha}}
{(2k-1-\alpha)!\cdot
\alpha!}\sum_{\gamma=0}^NC_{\gamma}(h\gamma)^{\alpha}, \eqno(25)
$$
$$ D=\frac{1}{4}\sum_{\gamma=0}^NC_{\gamma}e^{h\gamma}.
\eqno(26)
$$

We denote
$$
Q^{(-)}_{m-2}(h\beta)=P_{m-2}(h\beta)+Q_{m-2}(h\beta),\ \ a^-=d+D,
\eqno(27)
$$
$$
Q^{(+)}_{m-2}(h\beta)=P_{m-2}(h\beta)-Q_{m-2}(h\beta),\ \ a^+=d-D.
\eqno(28) $$ Then, taking into account (11), (16), (23), (27),
(28), we obtain following problem

\textbf{Problem B}. \emph{Find solution of the equation}
$$
D_m(h\beta)*u(h\beta)=0\mbox{  when } h\beta\not\in [0,1]
$$
\emph{which has view:}
$$
u(h\beta)= \left\{
\begin{array}{ll}
-\frac{e^{h\beta}}{4}q_{m-1}
+Q^{(2m-3)}(h\beta)+Q_{m-2}^{(-)}(h\beta)+a^-e^{-h\beta}, & \beta<0;\\
f_m(h\beta), & 0\leq \beta\leq N;\\
\frac{e^{h\beta}}{4}q_{m-1}-Q^{(2m-3)}(h\beta)+Q_{m-2}^{(+)}(h\beta)
+a^+e^{-h\beta}, & \beta>N.
\end{array}
\right. \eqno(29)
$$
Here polynomials $Q_{m-2}^{(-)}(h\beta),$ $Q_{m-2}^{(+)}(h\beta)$
and constants $a^-$, $a^+$ are unknowns. If we will find them,
then from (27), (28) we have
$$
P_{m-2}(h\beta)=\frac{1}{2}(Q_{m-2}^{(-)}(h\beta)+
Q_{m-2}^{(+)}(h\beta)),\ \ d=\frac{1}{2}[a^-+a^+],
$$
$$
Q_{m-2}(h\beta)=\frac{1}{2}(Q_{m-2}^{(-)}(h\beta)-
Q_{m-2}^{(+)}(h\beta)),\ \ D=\frac{1}{2}[a^- -a^+].
$$
Unknowns $Q_{m-2}^{(-)}(h\beta)$, $Q_{m-2}^{(+)}(h\beta)$, $a^-$
and $a^+$ are found from conditions (22). Thus, problem B and
respectively problem A are simultaneously solved.

We introduce following designations
$$
 a_k=\frac{A_k}{\lambda_kp}\sum\limits^\infty_{\gamma=1}
\lambda^\gamma_k\left[-\frac{e^{-h\gamma}}
{4}q_{m-1}+Q^{2m-3}(-h\gamma)+Q^{(-)}_{m-2}(-h\gamma)+a^-e^{h\gamma}-f_m(-h\gamma)\right],
\eqno (30) $$
$$
b_k=\frac{A_k}{\lambda_kp}\sum\limits^\infty_{\gamma=1}
\lambda^\gamma_k\bigg[-\frac{e^{h\gamma+1}}
{4}q_{m-1}-Q^{2m-3}(1+h\gamma)+Q^{(+)}_{m-2}(1+h\gamma)+
$$
$$
+a^+e^{-1-h\gamma}-f_m(1+h\gamma)\bigg], \eqno(31)
$$
$$
a=-\frac{2e^h}{p}\left[-\frac{e^{-h}}{4}q_{m-1}+Q^{(2m-3)}(-h)+Q_{m-2}^{(-)}(-h)+a^-e^h-f_m(-h)\right],
\eqno(32)
$$
$$
b=-\frac{2e^h}{p}\left[\frac{e^{1+h}}{4}q_{m-1}-Q^{(2m-3)}(1+h)+Q_{m-2}^{(+)}(1+h)+a^+
e^{-h-1}-f_m(1+h)\right]. \eqno(33)
$$

In present paper, using above mentioned algorithm [2], following
theorems are proved.

\textbf{Theorem 3.} \emph{The coefficients of the weight optimal
quadrature formulas of the form  (1) in the $W_2^{(m,m-1)}(0,1)$
space have view:}
$$
C_{\beta}= \left\{
\begin{array}{ll}
D_m(h\beta)*f_m(h\beta)|_{\beta=0}+a+\sum\limits_{k=1}^{m-1}(a_k+b_k\lambda_k^N),&\beta=0;\\
D_m(h\beta)*f_m(h\beta)+\sum\limits_{k=1}^{m-1}(a_k\lambda_k^{\beta}+b_k\lambda_k^{N-\beta}),
              &\beta=\overline{1,N-1};\\
D_m(h\beta)*f_m(h\beta)|_{\beta=N}+b+\sum\limits_{k=1}^{m-1}(a_k\lambda_k^N+b_k),&\beta=N.\\
\end{array}
\right. \eqno(34)
$$
\emph{where $D_m(h\beta)$ is defined by (17)}, $f_m(h\beta)$
\emph{is defined by} (15), \emph{$a_k$, $b_k$, $a$, $b$ are
determined by (30)-(33), $\lambda_k$ are roots of the polynomial
$P_{2m-2}(\lambda)$, which given by equality (19) and
$|\lambda_k|<1$.}

\textbf{Theorem 4.} \emph{The optimal coefficients of the
quadrature formulas of the form  (1) in the $W_2^{(1,0)}(0,1)$
space when $p(x)=1$ have following view:}
$$
C_{\beta}=\left\{
\begin{array}{l}
\frac{e^h-1}{e^h+1},\ \ \beta=0,N; \\
\frac{2(e^h-1)}{e^h+1},\ \ \beta=\overline{1,N-1},
\end{array}
\right.
$$
\emph{where} $h=\frac{1}{N}$, $N=1,2,...$.

\textbf{Theorem 5.} \emph{The optimal coefficients of the
quadrature formulas of the form (1) in the  $W_2^{(2,1)}(0,1)$
space when $p(x)=1$ have following view:}
$$
C_\beta=\left\{ \begin{array}{ll}
1-\frac{h}{e^h-1}-K(h)(\lambda_1-\lambda_1^N), & \beta=0\\
h+K(h)\left((e^h-\lambda_1)\lambda_1^\beta+(1-\lambda_1e^h)\lambda_1^{N-\beta}\right),
& \beta=\overline{1,N-1}\\
-1+e^h\left(\frac{h}{e^h-1}-K(h)(\lambda_1-\lambda_1^N)\right), &
\beta=N.
\end{array}
\right.
$$
\emph{where} $$
K(h)=\frac{(2e^h-2-he^h-h)(\lambda_1-1)}{2(e^h-1)^2(\lambda_1+\lambda_1^{N+1})}$$
$$
\lambda_1=\frac{h(e^{2h}+1)-e^{2h}+1-(e^h-1)\sqrt{h^2(e^h+1)^2+2h(1-e^h)}}
{1-e^{2h}+2he^h},\ \ |\lambda_1|<1,
$$
\emph{$h=\frac{1}{N}, N=1,2,....$ }\\[0.5cm]

\textbf{Proofs of theorems}

\textbf{Proof of theorem 3. } We obtained, that
$$
C_\beta=D_m(h\beta)\ast u(h\beta).
$$
Hence, taking into account (17), (29) and properties of the
operator $D_m(h\beta)$ we have
$$
C_\beta=D_m(h\beta)\ast
u(h\beta)=\sum\limits_{\gamma=-\infty}^{\infty}D_m(h\beta-h\gamma)u(h\gamma)=
$$
$$
=\sum\limits_{\gamma=1}^{\infty}D_m(h\beta+h\gamma)
\left[-\frac{e^{-h\gamma}}{4}q_{m-1}+Q^{(2m-3)}(-h\gamma)+Q_{m-2}^{(-)}(-h\gamma)+
a^-e^{h\gamma}\right]+
$$
$$
+\sum\limits_{\gamma=0}^ND_m(h\beta-h\gamma)f_m(h\gamma)+\sum\limits_{\gamma=1}^{\infty}
D_m(h(\gamma+N)-h\beta)\Bigg[\frac{e^{h\gamma+1}}{4}q_{m-1}-Q^{(2m-3)}(1+h\gamma)+
$$
$$
+Q_{m-2}^{(+)}(1+h\gamma)+a^+e^{-1-h\gamma}\Bigg].
$$

According to the definition of convolution of discrete functions
we get
$$
C_\beta=\sum\limits^\infty_{\gamma=-\infty}D_m(h\beta-h\gamma)f_m(h\gamma)+
\sum\limits^\infty_{\gamma=1}D_m(h\beta+h\gamma)\bigg[-\frac{e^{-h\gamma}}{4}q_{m-1}+
$$
$$
+Q^{(2m-3)}(-h\gamma)+Q_{m-2}^{(-)}(-h\gamma)+a^-e^{h\gamma}-f_m(-h\gamma)\bigg]+
\sum\limits_{\gamma=1}^{\infty}D_m(h(\gamma+N)-h\beta)\times
$$
$$
\times\bigg[\frac{e^{h\beta}}{4}q_{m-1}-Q^{(2m-3)}(1+h\gamma)+Q_{m-2}^{(+)}(h\gamma+1)+
a^+e^{-h\gamma-1}-f_m(h\gamma+1)\bigg]=
$$
$$
=D_m(h\beta)*f_m(h\beta)+\sum\limits_{\gamma=1}^{\infty}D_m(h\beta+h\gamma)\bigg[
-\frac{e^{-h\gamma}}{4}q_{m-1}+Q^{(2m-3)}(-h\gamma)+Q_{m-2}^{(-)}(-h\gamma)+
$$
$$
+a^-e^{h\gamma}-f_m(-h\gamma)\bigg]+\sum\limits_{\gamma=1}^{\infty}D_m(h(N+\gamma)-
h\beta)\bigg[\frac{e^{h\gamma+1}}{4}q_{m-1}-Q^{(2m-3)}(1+h\gamma)+
$$
$$
+Q_{m-2}^{(+)}(1+h\gamma)+a^+e^{-h\gamma}-f_m(1+h\gamma)\bigg].
\eqno(35)
$$
Hence when $\beta=0$ using equality (17) we obtain
$$
C_0=D_m(h\beta)\ast
f_m(h\beta)|_{\beta=0}+\frac{-2e^h}{p}\bigg[-\frac{e^{-h}}{4}q_{m-1}+Q^{(2m-3)}(-h)+
Q_{m-2}^{(-)}(-h)+
$$
$$
+a^-e^h-f_m(-h)\bigg]+\sum\limits_{k=1}^{m-1}\frac{A_k}{\lambda_kp}\sum\limits_{\gamma=1}
^{\infty}\lambda_k^{\gamma}\bigg[-\frac{e^{-h\gamma}}{4}q_{m-1}+Q^{(2m-3)}(-h\gamma)+
$$
$$
+Q_{m-2}^{(-)}(-h\gamma)+a^-e^{h\gamma}-f_m(-h\gamma)\bigg]+\sum\limits_{k=1}^{m-1}
\frac{A_k\lambda_k^N}{\lambda_kp}\sum\limits_{\gamma=1}^{\infty}\bigg[
\frac{e^{1+h\gamma}}{4}q_{m-1}-
$$
$$
-Q^{(2m-3)}(1+h\gamma)+Q^{(+)}(1+h\gamma)+a^+e^{-1-h\gamma}-f_m(1+h\gamma)\bigg].
\eqno(36) $$ From (36), taking account of (30), (31) and (32) we
have
$$
C_0=D_m(h\beta)*f_m(h\beta)|_{\beta=0}+a+\sum\limits_{k=1}^{m-1}(a_k+b_k\lambda_k^N).
$$

Let now $1\leq \beta\leq N-1$, then from (35) and (17), we get
$$
C_{\beta}=D_m(h\beta)*f_m(h\beta)+\sum\limits_{\gamma=1}^{\infty}\sum\limits_{k=1}^{m-1}
\frac{A_k}{p}\lambda_k^{\beta+\gamma-1}\bigg[-\frac{e^{-h\gamma}}{4}q_{m-1}+Q^{(2m-3)}(-h\gamma)+
$$
$$
+Q_{m-2}^{(-)}(-h\gamma)+a^-e^{h\gamma}-f_m(-h\gamma)\bigg]+\sum\limits_{\gamma=1}^{\infty}
\sum\limits_{k=1}^{m-1}\frac{A_k}{p}\lambda_k^{N+\gamma-\beta}\bigg[\frac{e^{h\gamma}}{4}q_{m-1}-
$$
$$
-Q^{(2m-3)}(1+h\gamma)+Q_{m-2}^{(+)}(1+h\gamma)+a^{(+)}e^{-1-h\gamma}-f_m(1+h\gamma)\bigg].
$$
Hence keeping in mind (30) and (31) we obtain
$$
C_\beta=D_m(h\beta)\ast
f_m(h\beta)+\sum\limits^{m-1}_{k=1}\left(a_k\lambda^\beta_k+b_k\lambda^{N-\beta}_k\right).
$$
Now consider case when  $\beta=N$. From (35), taking account of
equality (17) and designations (30), (31), (33), we get
$$
C_N=D_m(h\beta)\ast
f_m(h\beta)|_{\beta=N}+\sum\limits_{\gamma=1}^{\infty}
\sum\limits_{k=1}^{m-1}\frac{A_k}{\lambda_kp}\lambda_k^{N+\gamma}\bigg[-\frac{e^{-h\gamma}}
{4}q_{m-1}+Q^{(2m-3)}(-h\gamma)+
$$
$$
+Q_{m-2}^{(-)}(-h\gamma)+a^-e^{h\gamma}-f_m(-h\gamma)\bigg]+
\frac{-2e^h}{p}\bigg[\frac{e^{1+h}}{4}q_{m-1}-Q^{(2m-3)}(1+h)+
$$
$$
+Q_{m-2}^{(+)}(1+h)+a^+
e^{-h-1}-f_m(1+h)\bigg]+\sum\limits_{\gamma=1}^{\infty}\sum\limits_{k=1}^{m-1}
\frac{A_k}{\lambda_kp}\lambda_k^{\gamma}\bigg[-\frac{e^{1+h\gamma}}{4}q_{m-1}-
$$
$$
-Q^{(2m-3)}(1+h\gamma)+Q_{m-2}^{(+)}(1+h\gamma)+a^+e^{-1-h\gamma}-f_m(1+h\gamma)\bigg]=
$$
$$
=D_m(h\beta)*f_m(h\beta)|_{\beta=N}+b+\sum\limits_{k=1}^{m-1}(a_k\lambda_k^N+b_k).
$$
Thus, we get the statement of theorem. Theorem 3 is proved.

Consider case when $p(x)=1$. Then system (11)-(14) of Wiener -Hopf
type has following view $$
\left\{
\begin{array}{l}
\sum\limits_{\gamma=0}^NC_{\gamma}\psi_m(h\beta-h\gamma)+P_{m-2}(h\beta)+de^{-h\beta}=f_m(h\beta),
\ \ \beta=\overline{0,N};\\
\sum\limits_{\gamma=0}^NC_{\gamma}(h\gamma)^{\alpha}=\frac{1}{\alpha+1},\
\ \alpha=\overline{0,m-2};\\
\sum\limits_{\gamma=0}^NC_{\gamma}e^{-h\gamma}=1-e^{-1}.
\end{array}
\right. \eqno (37) $$ where $\psi_m(x)$ is defined by equality
(5), $P_{m-2}(h\beta)$ is polynomial of degree $m-2$ of $h\beta$,
$d$ is unknown coefficient
$$
f_m(h\beta)=\frac{e^{h\beta}+e^{-h\beta}+e^{1-h\beta}+e^{h\beta-1}-4}{4}-
\sum_{k=1}^{m-1}\frac{(h\beta)^{2k}+(1-h\beta)^{2k}}{2\cdot
(2k)!}. \eqno (38)
$$

\textbf{Proof of theorem 4.} Let $m=1$, then from (37) we have
$$
\sum\limits_{\gamma=0}^NC_{\gamma}\frac{\mathrm{sign}(h\beta-h\gamma)}{4}
(e^{h\beta-h\gamma}-e^{h\gamma-h\beta})+de^{-h\beta}=f_1(h\beta),\
\beta=\overline{0,N}, \eqno (39)
$$
$$
\sum\limits_{\gamma=0}^NC_{\gamma}e^{-h\gamma}=1-e^{-1}. \eqno
(40)
$$

According to (34) in this case the optimal coefficients have
following form
$$
C_{\beta}= \left\{
\begin{array}{ll}
D_1(h\beta)*f_1(h\beta)|_{\beta=0}+a,&\beta=0;\\
D_1(h\beta)*f_1(h\beta),&\beta=\overline{1,N-1};\\
D_1(h\beta)*f_1(h\beta)|_{\beta=N}+b,&\beta=N.\\
\end{array}
\right. \eqno(41)
$$

To find optimal coefficients (41) we need to calculate the
convolution $D_1(h\beta)\ast f_1(h\beta)$. For this using
equalities (17), (38) and theorem 2 we obtain
$$
D_1(h\beta)*f_1(h\beta)=D_1(h\beta)*\left(\frac{e^{h\beta}+e^{-h\beta}+e^{1-h\beta}+
e^{h\beta-1}}{4}-1\right)=
$$
$$
\frac{1}{4}\left(D_1(h\beta)*e^{h\beta}+D_1(h\beta)*e^{-h\beta}+
D_1(h\beta)*e^{1-h\beta}+D_1(h\beta)*e^{h\beta-1}\right)-D_1(h\beta)*1=
$$
$$
=-D_1(h\beta)*1=-\sum\limits_{\gamma=-\infty}^{\infty}D_1(h\gamma)=-2D_1(h)-D_1(0)=
$$
$$
=-2\cdot\frac{1}{p}(-2e^h)-\frac{1}{p}\cdot
2C=-\frac{1}{p}(-4e^h+2(1+e^{2h}))=\frac{2(e^h-1)}{e^h+1}.
$$
Then from (41) we have
$$
C_{\beta}= \left\{
\begin{array}{ll} \frac{2(e^h-1)}{e^h+1}+a, &
\beta=0\\
\frac{2(e^h-1)}{e^h+1}, & \beta=\overline{1,N-1}\\
\frac{2(e^h-1)}{e^h+1}+b, & \beta=N.
\end{array}
\right. \eqno(42)
$$
Hence one can see, that for finding of coefficients we need to
determine $a$ and $b$. For this directly using (42), (40), from
(39) taking account of $\beta=\overline{1,N-1}$, we have
$$
S=\sum\limits^N_{\gamma=0}C_\gamma\frac{\mathrm{sign}(h\beta-h\gamma)}{4}(e^{h\beta-h\gamma}
-e^{h\gamma-h\beta})=\frac{1}{2}\sum\limits^\beta_{\gamma=0}C_\gamma(e^{h\beta-h\gamma}
-e^{h\gamma-h\beta})-
$$
$$
-\frac{1}{4}\sum\limits^N_{\gamma=0}C_\gamma(e^{h\beta-h\gamma}
-e^{h\gamma-h\beta})=\frac{1}{2}\left(\sum\limits^\beta_{\gamma=0}c_\gamma
e^{h\beta-h\gamma} -\sum\limits^\beta_{\gamma=0}C_\gamma
e^{h\gamma-h\beta}\right)-
$$
$$
-\frac{1}{4}\left(e^{h\beta}(1-e^{-1})-e^{-h\beta}
\sum\limits^N_{\gamma=0}C_\gamma e^{h\gamma}\right)=
$$
$$
=\frac{1}{2}\left(a(e^{h\beta}-e^{-h\beta})+\sum\limits^\beta_{\gamma=0}
\frac{2(e^h-1)}{e^h+1}(e^{h\beta-h\gamma}
-e^{h\gamma-h\beta})\right)-
$$
$$
-\frac{1}{4}\left(e^h\beta(1-e^{-1})-e^{-h\beta}(a+be+\sum\limits^N_{\gamma=0}
\frac{2(e^h-1)}{e^h+1} e^{h\gamma})\right)=
$$
$$
=\frac{1}{2}\left(a(e^{h\beta}
-e^{-h\beta})+\frac{2(e^h-1)}{e^h+1}\left(e^{h\beta}\frac{1-(e^{-h})^{\beta+1}}{1-e^{-h}}
-e^{-h\beta}\frac{1-(e^h)^{\beta+1}}{1-e^h}\right)\right)-
$$
$$
-\frac{1}{4}\left(e^{h\beta}(1-e^{-1})-e^{-h\beta}(a+be+\frac{2(e^h-1)}{e^h+1}
\frac{1-(e^h)^{N+1}}{1-e^h})\right)=
$$
$$
=\frac{1}{2}\left(a(e^{h\beta}-e^{-h\beta})+\frac{2(e^h-1)}{e^h+1}\left(e^{h\beta}
\frac{e^h-e^{-h\beta}}{e^h-1}+e^{-h\beta}\frac{1-e^{h\beta+h}}{e^h-1}\right)\right)-
$$
$$
-\frac{1}{4}\left(e^{h\beta}(1-e^{-1})-e^{-h\beta}
\left(a+be+\frac{2(e^{h+1}-1)}{e^h+1}\right)\right)=
$$
$$
=e^{h\beta}\left(\frac{a}{2}+\frac{2e^h}{e^h+1}-\frac{1-e^{-1}}{4}\right)+e^{-h\beta}\left(
-\frac{a}{2}+\frac{1}{e^h+1}+\frac{a+be}{4}+\frac{e^{h+1}-1}{2(e^h+1)}\right)-1.
$$

Substituting obtained expression for $S$ into (39) we get
$$
e^{h\beta}\left(\frac{a}{2}+\frac{e^h}{e^h+1}-\frac{e-1}{4e}\right)+e^{-h\beta}
\left(\frac{be-a}{4}+\frac{e^{h+1}+1}{2(e^h+1)}\right)-1+de^{-h\beta}=
$$
$$
=e^{h\beta}(\frac{1+e^{-1}}{4})+e^{-h\beta}\frac{1-e}{4}-1.
$$
Hence and from (40) we obtain
$$
\left\{ \begin{array}{l}
\frac{a}{2}+\frac{e^h}{e^h+1}-\frac{e-1}{4e}=\frac{e+1}{4e}\\
\frac{be-a}{4}+\frac{e^{h+1}+1}{2(e^h+1)}+d=\frac{1+e}{4}\\
ae+b+\frac{2(e^{h+1}-1}{e^h+1}=e-1
\end{array}
\right.
$$
Solving this system, we will find unknowns $a, b, d,$ i.e.$$
a=\frac{1-e^h}{e^h+1}, \ b=\frac{1-e^h}{e^h+1}, \ d=0. \eqno
(43)$$ Thus, substituting (43) into (42), we obtain the assertion
of theorem. Theorem 4 is proved.

\textbf{Proof of theorem 5.} In $m=2$ from (34) we obtain
$$
C_\beta=\left\{
\begin{array}{ll} D_2(h\beta)\ast
f_2(h\beta)|_{\beta=0}+a+a_1+b_1\lambda_1^N, & \beta=0\\
D_2(h\beta)\ast
f_2(h\beta)+a_1\lambda_1^\beta+b_1\lambda_1^{N-\beta}, &
1\leq\beta\leq N-1\\
D_2(h\beta)\ast f_2(h\beta)|_{\beta=N}+b+a_1\lambda_1^N+b_1, &
\beta=N\\
\end{array}
\right. \eqno (44)
$$

First we will calculate the convolution $D_2(h\beta)\ast
f_2(h\beta)$. Using equalities  (17), (38) and theorem 2 we get
$$
D_2(h\beta)\ast
f_2(h\beta)=D_2(h\beta)\ast\bigg(\frac{e^{h\beta}+e^{-h\beta}+e^{1-h\beta}+e^{h\beta-1}-4}
{4}-
$$
$$ -\frac{(h\beta)^2+(1-h\beta)^2}{4}\bigg)=
-D_2(h\beta)\ast\left(\frac{2(h\beta)^2}{4}\right)=-\frac{1}{2}D_2(h\beta)\ast(h\beta)^2.
\eqno(45)
$$ From equality (20) we obtain, that polynomial
$$
P_2(\lambda)=\lambda^2(1-e^{2h}+2he^h)+(2(e^{2h}-1)-2h(e^{2h}+1))\lambda+(1-e^{2h}+2he^h)
$$
has two roots  $\lambda_1$ and $\lambda_2$, $\lambda_1\lambda_2=1$
and $|\lambda_1|<1$
$$
\lambda_1+\lambda_2=\frac{2h(e^{2h}+1)-2(e^{2h}-1)}{(1-e^{2h}+2he^h)},
\ p=1-e^{2h}+2he^h \eqno (46)
$$
Further,
$$
D_2(h\beta)\ast
(h\beta)^2=\sum\limits^\infty_{\gamma=-\infty}D_2(h\gamma)(h\beta-h\gamma)^2=
\sum\limits^\infty_{\gamma=-\infty}D_2(h\gamma)(h\gamma)^2=
$$
$$
2D_2(h)h^2+2\sum\limits^\infty_{\gamma=2}D_2(h\gamma)(h\gamma)^2=
2h^2(\frac{1}{p}(-2e^h+A_1))+
$$
$$
+2h^2\sum\limits^\infty_{\gamma=2}\frac{A_1}{p\lambda_1}\lambda_1^\gamma\gamma^2=
\frac{2h^2}{p}\left[-2e^h+\frac{A_1}{\lambda_1}\sum\limits^\infty_{\gamma=1}\lambda_1^\gamma
\gamma^2\right]=
$$
$$
=\frac{2h^2}{p}\left[-2e^h+\frac{A_1}{\lambda_1}\frac{1}{1-\lambda_1}\sum\limits^2_{i=0}
(\frac{\lambda_1}{1-\lambda_1})^i\Delta^i0^2\right]=
$$
$$
=\frac{2h^2}{p}\left[-2e^h+\frac{A_1}{\lambda_1}\frac{1}{1-\lambda_1}\left(
\frac{\lambda_1}{1-\lambda_1}+\frac{2\lambda_1^2}{(1-\lambda_1)^2}\right)\right]=
$$
$$
=\frac{2h^2}{p}\left[-2e^h+\frac{A_1}{\lambda_1}\frac{1}{1-\lambda_1}
\frac{\lambda_1+\lambda_1^2}{(1-\lambda_1)^2}\right]=
$$
$$
=\frac{2h^2}{p}\left[-2e^h+\frac{A_1}{\lambda_1}\frac{\lambda_1(1+\lambda_1)}{(1-\lambda_1)^3}
\right]=\frac{2h^2}{p}\left[-2e^h+\frac{A_1(1+\lambda_1)}{(1-\lambda_1)^3}\right].
$$
Hence, taking into account equalities (19) and (46), we have
$$
D_2(h\beta)\ast
(h\beta)^2=\frac{2h^2}{p}\left[-2e^h-\frac{2(\lambda_1)^2(\lambda_1(e^{2h}+1)-e^h(\lambda_1^2+1))(1+\lambda_1)}
{(\lambda_1^2-1)(1-\lambda_1)^3}\right]=
$$
$$
\frac{2h^2}{p}\left[-2e^h-\frac{2(\lambda_1(e^{2h}+1)-e^h(\lambda_1^2+1))}{(1-\lambda_1)^2}\right]=
$$
$$
=-\frac{4h^2}{p}\left[e^h+\frac{\lambda_1(e^{2h}+1)-e^h(\lambda_1^2+1)}{(1-\lambda_1)^2}\right]=
$$
$$
=-\frac{4h^2}{p}\left[\frac{e^h-2\lambda_1e^h+e^h\lambda_1^2+\lambda_1(e^{2h}+1)-e^h(\lambda_1^2+1)}
{(1-\lambda_1)^2}\right]=
$$
$$
=-\frac{4h^2}{p}\frac{(e^h-1)^2\lambda_1}{(1-\lambda_1)^2}=\frac{-4h^2(e^h-1)^2}{p}
\frac{1}{\lambda_2(1-2\lambda_1+\lambda_1^2)}=
$$
$$
=-\frac{4h^2(e^h-1)^2}{p}\frac{1}{\lambda_2-2+\lambda_1}=-\frac{4h^2(e^h-1)^2}{p}
\frac{1}{\frac{2h(e^{2h}+1)-2(e^{2h}-1)}{1-e^{2h}+2he^h}-2}=
$$
$$
=-\frac{4h^2(e^h-1)^2}{p}\frac{p}{2h(e^{2h}+1)-2(e^{2h}-1)-2(1-2e^{2h})-4he^h}=
$$
$$
=-\frac{4h^2(e^h-1)^2}{2h(e^{2h}-2e^h+1)}=-\frac{4h^2(e^h-1)^2}{2h(e^h-1)^2}=-2h.
$$
So $$ D_2(h\beta)\ast(h\beta)^2=-2h. \eqno (47) $$ If we
substitute equality (47)  into (45) we obtain
$$
D_2(h\beta)\ast f_2(h\beta)=-\frac{1}{2}D_2(h\beta)\ast
(h\beta)^2=-\frac{1}{2}(-2h)=h. \eqno(48) $$ Taking into account
(48) from (44) we get
$$
 C_\beta=\left\{
\begin{array}{ll} h+a+a_1+b_1\lambda_1^N,
&\beta=0,\\
h+a_1\lambda_1^\beta+b_1\lambda_1^{N-\beta}, & 1\leq\beta\leq
N-1,\\
h+b+a_1\lambda_1^N+b_1, & \beta=N.
\end{array}
\right. \eqno (49)
$$
Here $a_1,b_1,a,b$ are unknown constants. If we find these
constants, then we obtain explicit form of optimal coefficients.

When $m=2$ the system (8)-(10) of Wiener-Hopf type for optimal
coefficients have following view
$$
\sum\limits_{\gamma=0}^N\frac{\mathrm{sign}(h\beta-h\gamma)}{2}\left(
\frac{e^{h\beta-h\gamma}-e^{h\gamma-h\beta}}{2}-(h\beta-h\gamma)\right)+de^{-h\beta}=
$$
$$
=\frac{e^{h\beta}+e^{-h\beta}+e^{1-h\beta}+e^{h\beta-1}-4}{4}-\frac{(h\beta)^2+(1-h\beta)^2}{4},
\eqno(50)
$$
$$
\sum\limits_{\gamma=0}^NC_{\gamma}=1,\eqno(51)$$
$$
\sum\limits_{\gamma=0}^NC_{\gamma}e^{-h\gamma}=1-e^{-1}. \eqno(52)
$$
To find explicit form of optimal coefficients, substituting
equality (49) into (50), we obtain identity by variable $h\beta$.
Equating to zero coefficients of $e^{h\beta}$, $e^{-h\beta}$,
$h\beta$, $(h\beta)^0$, we get system for unknowns $a_1,b_1$ and
solving it we will find $a_1,b_1$.

First from orthogonality conditions (51), (52) coefficients $C_0$
and $C_N$ we will express by $C_{\gamma}$,
$\gamma=\overline{1,N-1}$
$$
\left\{
\begin{array}{l}
C_0+C_N=1-\sum\limits_{\gamma=1}^{N-1}C_{\gamma},\\
C_0+C_Ne^{-1}=1-e^{-1}-\sum\limits_{\gamma=1}^{N-1}C_{\gamma}e^{-h\gamma}.
\end{array}
\right.
$$
Hence
$$
C_0=\frac{e-2}{e-1}+\sum\limits_{\gamma=1}^{N-1}C_{\gamma}\frac{1-e^{1-h\gamma}}{e-1},
\eqno(53)
$$
$$
C_N=\frac{1}{e-1}+\sum\limits_{\gamma=1}^{N-1}C_{\gamma}\frac{e^{1-h\gamma}-e}{e-1}.
\eqno(54)
$$
Now consider first sum of the equality (50) and taking into
account definition of $\mathrm{sign}(x)$ we get
$$
g(h\beta)=\sum\limits_{\gamma=0}^NC_{\gamma}\frac{\mathrm{sign}(h\beta-h\gamma)}{2}
\left(\frac{e^{h\beta-h\gamma}-e^{h\gamma-h\beta}}{2}-(h\beta-h\gamma)\right)=
$$
$$
=\sum\limits_{\gamma=0}^{\beta}C_{\gamma}\left(\frac{e^{h\beta-h\gamma}-e^{h\gamma-h\beta}}{2}-
(h\beta-h\gamma)\right)-
$$
$$
-\frac{1}{2}\sum\limits_{\gamma=0}^NC_{\gamma}
\frac{\mathrm{sign}(h\beta-h\gamma)}{2}
\left(\frac{e^{h\beta-h\gamma}-e^{h\gamma-h\beta}}{2}-(h\beta-h\gamma)\right)=
$$
$$
=C_0\left(\frac{e^{h\beta}-e^{-h\beta}}{2}-h\beta\right)+\sum\limits_{\gamma=1}^{\beta}
C_{\gamma}\left(\frac{e^{h\beta-h\gamma}-e^{h\gamma-h\beta}}{2}-(h\beta-h\gamma)\right)-
$$
$$
-\frac{1}{2}\left(\frac{e^{h\beta}}{2}\sum\limits_{\gamma=0}^NC_{\gamma}e^{-h\gamma}-
\frac{e^{-h\beta}}{2}\sum\limits_{\gamma=0}^NC_{\gamma}e^{h\gamma}-
\sum\limits_{\gamma=0}^NC_{\gamma}h\beta+\sum\limits_{\gamma=0}^NC_{\gamma}h\gamma\right)=
$$
$$
=C_0\left(\frac{e^{h\beta}-e^{-h\beta}}{2}-h\beta\right)+\sum\limits_{\gamma=1}^{\beta}
C_{\gamma}\left(\frac{e^{h\beta-h\gamma}-e^{h\gamma-h\beta}}{2}-(h\beta-h\gamma)\right)-
$$
$$
-\frac{1}{2}\left(\frac{e^{h\beta}}{2}(1-e^{-1})-
\frac{e^{-h\beta}}{2}\sum\limits_{\gamma=0}^NC_{\gamma}e^{h\gamma}-
h\beta+\sum\limits_{\gamma=0}^NC_{\gamma}h\gamma\right). \eqno
(55)
$$

Using (49) - form of optimal coefficients $C_{\gamma}$,
$\gamma=\overline{1,N-1}$, we calculate following sums
$$
S_1=\sum\limits_{\gamma=1}^{N-1}C_{\gamma}\frac{1-e^{1-h\gamma}}{e-1}=\frac{1}{e-1}\left(1-h
+a_1\frac{\lambda_1-\lambda_1^N}{1-\lambda_1}+b_1\frac{\lambda_1-\lambda_1^N}{1-\lambda_1}\right)-
$$
$$
-\frac{e}{e-1}\left(h\frac{1-e^{h-1}}{e^h-1}+a_1\frac{\lambda_1-\lambda_1^Ne^{h-1}}{e^h-
\lambda_1}+b_1\frac{\lambda_1^N-\lambda_1e^{h-1}}{\lambda_1e^h-1}\right),
$$
$$
S_2=\sum\limits_{\gamma=1}^{N-1}C_{\gamma}\frac{e^{1-h\gamma}-e}{e-1}=\frac{e}{e-1}
\Bigg(h\frac{1-e^{h-1}}{e^h-1}+a_1\frac{\lambda_1-\lambda_1^Ne^{h-1}}{e^h-\lambda_1}+
$$
$$
+b_1\frac{\lambda_1^N-\lambda_1e^{h-1}}{\lambda_1e^h-1}\Bigg)-\frac{e}{e-1}\left(
1-h+a_1\frac{\lambda_1-\lambda_1^N}{1-\lambda_1}+b_1\frac{\lambda_1^N-\lambda_1}{\lambda_1-1}
\right),
$$
$$
S_3=\sum\limits_{\gamma=1}^{\beta}C_{\gamma}\left(\frac{e^{h\beta-h\gamma}-
e^{h\gamma-h\beta}}{2}-(h\beta-h\gamma)\right)=
$$
$$
=\sum\limits_{\gamma=1}^{\beta}(h+a_1\lambda_1^{\gamma}+b_1\lambda_1^{N-\gamma})
\left(\frac{e^{h\beta-h\gamma}-
e^{h\gamma-h\beta}}{2}-(h\beta-h\gamma)\right)=
$$
$$
=\frac{e^{h\beta}}{2}\sum\limits_{\gamma=1}^{\beta}(h+a_1\lambda_1^\gamma+b_1\lambda_1^{N-\gamma})
e^{-h\gamma}-\frac{e^{-h\beta}}{2}\sum\limits_{\gamma=1}^\beta(h+a_1\lambda_1^\gamma+b_1\lambda_1^{N-\gamma})
e^{h\gamma}-
$$
$$
-h\beta\sum\limits_{\gamma=1}^\beta(h+a_1\lambda_1^\gamma+b_1\lambda_1^{N-\gamma})+
\sum\limits_{\gamma=1}^\beta(h+a_1\lambda_1^\gamma+b_1\lambda_1^{N-\gamma})h\gamma.
$$
Hence using properties of geometric progression we get
$$
S_3=\frac{e^{h\beta}}{2}\left(h\frac{1-e^{-h\beta}}{e^h-1}+a_1\frac{\lambda_1-\lambda_1^{\beta+1}
e^{-h\beta}}{e^h-\lambda_1}+b_1\frac{\lambda_1^N-\lambda_1^{N-\beta}e^{-h\beta}}
{\lambda_1e^h-1}\right)-
$$
$$
-\frac{e^{-h\beta}}{2}\left(h\frac{e^h-e^{h+h\beta}}{1-e^h}+a_1\frac{\lambda_1e^h-\lambda_1^{\beta+1}
e^{h\beta+h}}{1-\lambda_1e^h}+b_1\frac{\lambda_1^Ne^h-\lambda_1^{N-\beta}e^{h+h\beta}}
{\lambda_1 -e^h} \right)-
$$
$$
-h\beta\left(h\beta+a_1\frac{\lambda_1-\lambda_1^{\beta+1}}{1-\lambda_1}+
b_1\frac{\lambda_1^N-\lambda_1^{N-\beta}}{\lambda_1-1}\right)+h^2\frac{(\beta+1)\beta}{2}+
$$
$$
+a_1h\frac{\beta(\lambda_1^{\beta+2}-\lambda_1^{\beta+1})+\lambda_1-\lambda_1^{\beta+1}}
{(1-\lambda_1)^2}+b_1h\frac{\lambda_1^{N+1}-\lambda_1^{N+1-\beta}+\beta(\lambda_1^{N-\beta}
-\lambda_1^{N+1-\beta})}{(\lambda_1-1)^2}=
$$
$$
=\frac{e^h\beta}{2}\left(\frac{h}{e^h-1}+\frac{a_1\lambda_1}{e^h-\lambda_1}+
\frac{b_1\lambda_1^N}{\lambda_1e^h-1}\right)-\frac{e^{-h\beta}}{2}\bigg(\frac{he^h}
{1-e^h}+
$$
$$
+a_1\frac{\lambda_1e^h}{1-\lambda_1e^h}+b_1\frac{\lambda_1^Ne^h}{\lambda_1-e^h}\bigg)+
\frac{1}{2}\bigg(\frac{h}{1-e^h}+\frac{a_1\lambda_1^{\beta+1}}{\lambda_1-e^h}+
$$
$$
+\frac{b_1\lambda_1^{N-\beta}}{1-\lambda_1e^h}\bigg)+\frac{1}{2}\left(\frac{he^h}{1-e^h}
+\frac{a_1\lambda_1^{\beta+1}e^h}{1-\lambda_1e^h}+\frac{b_1\lambda_1^{N-\beta}e^h}
{\lambda_1-e^h}\right)-
$$
$$
-h\beta\left(h\beta+a_1\frac{\lambda_1-\lambda_1^{\beta+1}}{1-\lambda_1}+b_1\frac{\lambda_1^N-
\lambda_1^{N-\beta}}{\lambda_1-1}\right)+h^2\frac{(\beta+1)\beta}{2}+
$$
$$
+a_1h\frac{\beta(\lambda_1^{\beta+2}-\lambda_1^{\beta+1})+\lambda_1-\lambda_1^{\beta+1}}
{(1-\lambda_1)^2}+b_1h\frac{\lambda_1^{N+1}-\lambda_1^{N+1-\beta}+\beta(\lambda_1^{N-\beta}-
\lambda_1^{N+1-\beta})}{(\lambda_1-1)2}.
$$
$$
S_4=\sum\limits_{\gamma=1}^{N-1} C_\gamma
e^{h\gamma}=\sum\limits_{\gamma=1}^{N-1}
(h+a_1\lambda_1^\gamma+b_1\lambda_1^{N-\gamma})e^{h\gamma}=
$$
$$
=h\frac{e^h-e}{1-e^h}+a_1\frac{\lambda_1e^h-\lambda_1^Ne}{1-\lambda_1e^h}
+b_1\frac{\lambda_1^{N}e^h-\lambda_1e}{\lambda_1-e^h}.
$$
$$
S_5=h\sum\limits_{\gamma=1}^{N-1}C_\gamma\gamma=h\sum\limits_{\gamma=1}^{N-1}
(h+a_1\lambda_1^\gamma+b_1\lambda_1^{N-\gamma})\gamma=\frac{h(N-1)}{2}+
$$
$$
+a_1h\frac{\lambda_1-\lambda_1^{N+1}-N\lambda_1^N(1-\lambda_1)}
{(\lambda_1-1)^2}+b_1h\frac{\lambda_1^{N+1}-\lambda_1-N(\lambda_1^2-\lambda_1)}
{(\lambda_1-1)^2}.
$$
Now by using equalities (53), (54), (55) and sums $S_1, S_2, S_3,
S_4, S_5$, equating coefficients in front of $e^{h\beta}$ and
$h\beta$ in left and right sides of equality  (50), respectively,
after some calculations we obtain following linear system for
$a_1$ and $b_1$.
$$
\left\{ \begin{array}{l}
a_1\frac{(e^h-1)(\lambda_1^{N+1}-\lambda_1e)}
{(1-\lambda_1)(e-1)(e^h-\lambda_1)}+b_1\frac{(e^h-1)(\lambda_1-\lambda_1^{N+1}e)
}{(e-1)(\lambda_1-1)(\lambda_1e^h-1)}=\frac{2-h}{2}-\frac{h}{e^h-1},\\
a_1\frac{1}{e^h-\lambda_1}+b_1\frac{1}{\lambda_1e^h-1}=0.\\
\end{array}
\right.
$$
Hence, solving system, we have
$$
a_1=\frac{(2e^h-2-he^h-h)(\lambda_1-1)}{2(e^h-1)^2(\lambda_1+\lambda_1^{N+1})}(e^h-\lambda_1),
\eqno (56)
$$
$$
b_1=\frac{(2e^h-2-he^h-h)(\lambda_1-1)}{2(e^h-1)^2(\lambda_1+\lambda_1^{N+1})}(1-e^h\lambda_1).
\eqno (57)
$$
From (53), taking account of $S_1$ and $a_1, b_1$ for $C_0$ we
get
$$
C_0=\frac{e^h-1-h}{e^h-1}-
\frac{(2e^h-2-he^h-h)(\lambda_1^2+\lambda_1^N-\lambda_1-\lambda_1^{N+1})}
{2(e^h-1)^2(\lambda_1+\lambda_1^{N+1})}.\eqno (58)
$$
And  from equality  (54), taking into account $S_2$, (56), (57) we
obtain optimal coefficient $C_N$
$$
C_N=\frac{he^h-e^h+1}{e^h-1}-
\frac{e^h(2e^h-2-he^h-h)(\lambda_1^2+\lambda_1^N-\lambda_1-\lambda_1^{N+1})}
{2(e^h-1)^2(\lambda_1+\lambda_1^{N+1})}.\eqno (59)
$$
We denote by
$$
K(h)=\frac{(2e^h-2-he^h-h)(\lambda_1-1)}{2(e^h-1)^2(\lambda_1+\lambda_1^{N+1})}
$$
and taking into account equalities (56), (57), (58) and (59) we
obtain assertion of theorem. Theorem 5 is proved.\\[0.5cm]

\textbf{References}
\begin{enumerate}
\item Shadimetov Kh.M., Hayotov A.R. Weight optimal quadrature formulas
in $W_2^{(m,m-1)}(0,1)$ space. Uzbek Mathematical Journal, 2002,
\No.3-4. pp.92-103. (in Russian)
\item
Sobolev S.L. Introduction to the Theory of Cubature Formulas. M.:
Nauka, 1974. - 808 p.
\item Shadimetov Kh.M., Hayotov A.R. Construction of Discrete
Analog of the Differential Operator
$\frac{d^{2m}}{dx^{2m}}-\frac{d^{2m-2}}{dx^{2m-2}}$. Uzbek
Mathematical Journal, 2004. \No 2 pp.85-95. (in Russian)
\item Shadimetov Kh.M., Hayotov A.R. Properties of the Discrete
Analogue of the Differential operator
$\frac{d^{2m}}{dx^{2m}}-\frac{d^{2m-2}}{dx^{2m-2}}$.
Kh.M.Shadimetov, A.R.Hayotov, 2004. pp.72-83. (in Russian)
\item
Sobolev S.L., Vaskevich V.L. Cubature Formulas. -Novosibirsk:
Institute of Mathematics SB of RAS, 1996, -484 p. (in Russian)
\end{enumerate}
\vspace{0.5cm}

\noindent
Kholmat Makhkambaevich Shadimetov\\
Institute of Mathematics and Information Technologies\\
Uzbek Academy of Sciences\\
Tashkent, 100125\\
Uzbekistan\\ [0.2cm]

\noindent
Abdullo Rakhmonovich Hayotov\\
Institute of Mathematics and Information Technologies\\
Uzbek Academy of Sciences\\
Tashkent, 100125\\
Uzbekistan\\
\textit{E-mail:} abdullo\_hayotov@mail.ru, hayotov@mail.ru.

\end{document}